\documentclass[11pt]{article}
\usepackage[utf8]{inputenc}
\usepackage[margin=1.2in]{geometry}

\usepackage{amsmath}
\usepackage{amssymb}
\usepackage{amsthm}
\usepackage{enumerate}
\usepackage{hyperref}
\usepackage{tikz-cd}

\DeclareMathOperator{\Hom}{Hom}
\DeclareMathOperator{\Tor}{Tor}
\DeclareMathOperator{\Ext}{Ext}
\DeclareMathOperator{\spec}{Spec}
\DeclareMathOperator{\maxspec}{mSpec}
\DeclareMathOperator{\coker}{coker}
\DeclareMathOperator{\im}{im}
\newcommand*{\mf}[1]{\mathfrak{#1}}
\newcommand*{\RMod}{R\text{-}\mathbf{Mod}}
\newcommand*{\To}{\Rightarrow}
\newcommand*{\NN}{\mathbb{N}}
\newcommand*{\ZZ}{\mathbb{Z}}
\newcommand*{\QQ}{\mathbb{Q}}

\newtheorem{theorem}{Theorem}[section]
\newtheorem{lemma}[theorem]{Lemma}
\newtheorem{corollary}[theorem]{Corollary}
\newtheorem{proposition}[theorem]{Proposition}
\newtheorem{fact}[theorem]{Fact}

\theoremstyle{definition}

\newtheorem{example}[theorem]{Example}

\theoremstyle{remark}
\newtheorem{remark}[theorem]{Remark}

\begin{document}

\title{Duality between injective envelopes and flat covers}

\author{Ville Puuska\thanks{Tampere University, Tampere, Finland. Email: puuskaville@gmail.com}}
\date{}

\maketitle

\begin{abstract}
We establish a duality between injective envelopes and flat covers over a commutative Noetherian ring. One case of this duality states that a morphism is an injective envelope, if and only if its Matlis dual is a flat cover. We also show that if we swap injective envelopes and flat covers in this duality, neither implication is true in general.
\end{abstract}

\section{Introduction}\label{sec: Introduction}

In this paper we establish the following duality between injective envelopes and flat covers of modules over a commutative Noetherian ring:
\begin{theorem}\label{thm: duality of inj envelope and flat cover}
Let \(R\) be a commutative Noetherian ring. Let \(i \colon M \to I\) be a morphism of \(R\)-modules. Then the following are equivalent:
\begin{enumerate}
\item \(i\) is an injective envelope;
\item the morphism \(\Hom_R(i,E)\) is a flat cover for some injective cogenerator \(E\);
\item the morphism \(\Hom_R(i,E)\) is a flat cover for all injective modules \(E\).
\end{enumerate}
\end{theorem}

Previously only partial dualities of this kind have been established. First, under the extra assumption that \(R\) is local, Belshoff and Xu showed in \cite[Theorem 4.1]{belshoff92} that \(R\) is complete and \(\dim R \leq 1\), if and only if all reflexive modules have reflexive flat covers and reflexive injective envelopes. They noted that this implies Theorem \ref{thm: duality of inj envelope and flat cover} whenever \(M\) is reflexive. This result was then generalized further in \cite[Theorem 11]{weimin96} and \cite[Corollary 3.4]{mao15}.

Later, Enochs and Huang showed in \cite[Theorem 3.7]{enochs12} that a (not necessarily commutative) ring is left Noetherian, if and only if a monomorphism \(i\) is an injective preenvelope precisely when \(i^+\) is a flat precover; here \((-)^+ := \Hom_\ZZ(-,\QQ/\ZZ)\). Combining this with Theorem \ref{thm: duality of inj envelope and flat cover}, we can characterize commutative Noetherian rings precisely as the commutative rings for which a morphism \(i\) is an injective envelope, if and only if \(i^+\) is a flat cover.

The key idea in the proof of Theorem \ref{thm: duality of inj envelope and flat cover} is to use the following alternative characterizations of injective envelopes and flat covers over a commutative Noetherian ring. It is well known that a monomorphism \(i \colon M \to E\), where \(E\) is injective, is an injective envelope, if and only if \(\Hom_{R_\mf p}(k(\mf p),i_\mf p)\) is an isomorphism for all \(\mf p \in \spec R\). Dailey proved in his thesis \cite[Proposition 4.2.7]{dailey16} a dual characterization of flat covers of cotorsion modules using the functor \(k(\mf p) \otimes_{R_\mf p} \Hom_R(R_\mf p,-)\). As our proof demonstrates, utilizing these two functors makes it far easier to attack the problem of proving the full duality between injective envelopes and flat covers.

We also consider the converse, i.e.~whether we can swap flat covers and injective envelopes in this duality. It turns out that for a wide variety of Noetherian rings this does not hold without extra restrictions. For example, we show that if \(R\) is a commutative Noetherian ring with \(\dim R \geq 2\), then there exists a flat cover \(f\) such that \(\Hom_R(f,E)\) is not an injective envelope for any injective cogenerator \(E\). However, the converse does hold with an added reflexivity requirement as was already shown in \cite[Corollary 3.3(2)]{mao15}.

\section{Preliminaries}\label{sec: Preliminaries}

Throughout, \(R\) is a commutative Noetherian ring with identity and modules are over \(R\) unless stated otherwise. We denote the category of \(R\)-modules by \(\RMod\). For any \(\mf p \in \spec R\), we denote \(k(\mf p) := R_\mf p/\mf pR_\mf p\).

Following \cite{enochs81}, we call a morphism \(f \colon F \to M\) of \(R\)-modules with \(F\) flat a \emph{flat precover}, if \(\Hom_R(F',F) \to \Hom_R(F',M) \to 0\) is exact for all flat modules \(F'\). If a flat precover \(f\) further satisfies \(fh = f\) only for automorphisms \(h \colon F \to F\), then we call \(f\) a \emph{flat cover}. Flat covers are unique up to isomorphism. Every module has a flat cover by \cite{bican01}.

For a module \(M\), the \emph{minimal flat resolution of \(M\)} is a resolution \(\cdots \xrightarrow{d_2} F_1 \xrightarrow{d_1} F_0 \xrightarrow{d_0} M \to 0\) such that \(F_i \to \im d_i\) is a flat cover for each \(i \in \NN\).

A module \(M\) is \emph{cotorsion}, if \(\Ext^1_R(F,M) = 0\) for all flat modules \(F\). Importantly, for any module \(M\) and injective module \(E\), \(\Hom_R(M,E)\) is pure injective and thus cotorsion \cite[Lemma 2.1]{enochs84}.

We will denote the Matlis dual of a module \(M\) by \(M^\vee := \Hom_R(M,\bigoplus_{\mf n} E(R/\mf n))\) where the direct sum is taken over the maximal ideals, i.e.~\(\mf n \in \maxspec R\). Similarly, for a morphism \(f \colon M \to N\), we set \(f^\vee := \Hom_R(f, \bigoplus_\mf n E(R/\mf n))\).

We will say that the module \(M\) is \emph{reflexive}, if the evaluation morphism \(M \to (M^\vee)^\vee\) is an isomorphism. Note that if \(M\) is reflexive, then all submodules and quotients of \(M\) are also reflexive.

To prove Theorem \ref{thm: duality of inj envelope and flat cover}, we need the following alternative characterizations of injective envelopes and flat covers. The characterization of injective envelopes is well known, but the characterization of flat covers was only recently proven by Dailey in his thesis.

\begin{fact}\label{fact: injective envelope iff soc isom}
Let \(i \colon M \to E\) be a morphism. Then, \(i\) is the injective envelope, if and only if \(E\) is injective, \(i\) is a monomorphism, and \(\Hom_{R_\mf p}(k(\mf p),i_\mf p)\) is an isomorphism for all \(\mf p \in \spec R\).
\end{fact}
\begin{proof}
The implication \(\To\) follows by using left exactness of \(\Hom_{R_\mf p}(k(\mf p),(-)_\mf p)\) and the fact that if \(0 \to M \to I^\bullet\) is a minimal injective resolution, then the morphisms in the complex \(\Hom_{R_\mf p}(k(\mf p),I^\bullet_\mf p)\) are all \(0\) (see e.g.~the proof of \cite[Proposition 3.2.9]{bruns98}).

For the implication \(\Leftarrow\) assume that \(E\) is injective, \(i\) is a monomorphism, and \(\Hom_{R_\mf p}(k(\mf p), i_\mf p)\) is an isomorphism for all \(\mf p \in \spec R\). Note that \(i\) factors through the injective envelope, i.e.~\(i \colon M \to E(M) \to E\), and the monomorphism \(E(M) \to E\) splits. If \(E(M) \to E\) is not an isomorphism, then we can split a factor \(E(R/\mf p) \subseteq E\) such that \(\im i \cap E(R/\mf p) = 0\). Since \(\Hom_{R_\mf p}(k(\mf p), E(R/\mf p)_\mf p) = k(\mf p) \neq 0\), this implies that \(\Hom_{R_\mf p}(k(\mf p), i_\mf p)\) is not surjective, which is of course a contradiction. Thus \(i\) is the injective envelope.
\end{proof}

\begin{fact}[{\cite[Proposition 4.2.7]{dailey16}}]\label{fact: flat cover iff top isom}
Let \(p \colon F \to M\) be an epimorphism where \(F\) is flat and cotorsion. Then, \(p\) is the flat cover, if and only if \(\ker p\) is cotorsion, and \(k(\mf p) \otimes_{R_\mf p} \Hom_R(R_\mf p,p)\) is an isomorphism for all \(\mf p \in \spec R\).
\end{fact}

\begin{lemma}\label{lem: duality of soc and top}
For \(\mf p \in \spec R\), an injective module \(E\), and a module \(M\), there exists a natural isomorphism
\[k(\mf p) \otimes_{R_\mf p} \Hom_R(R_\mf p,\Hom_R(M,E)) \cong \Hom_R(\Hom_{R_\mf p}(k(\mf p),M_\mf p),E).\]
\end{lemma}
\begin{proof}
The tensor-hom adjunction gives the natural isomorphism
\[k(\mf p) \otimes_{R_\mf p} \Hom_R(R_\mf p,\Hom_R(M,E)) \cong k(\mf p) \otimes_{R_\mf p} \Hom_R(M_\mf p,E).\]
Since \(k(\mf p)\) is finitely presented as an \(R_\mf p\)-module and \(E\) is injective, we have the natural isomorphism (see e.g.~\cite[Theorem 3.2.11]{enochs11})
\[k(\mf p) \otimes_{R_\mf p} \Hom_R(M_\mf p,E) \cong \Hom_R(\Hom_{R_\mf p}(k(\mf p),M_\mf p),E).\]
\end{proof}

\section{Duality between injective envelopes and flat covers}\label{sec: Duality between injective envelopes and flat covers}

We are now ready to prove our main result. Note that \(2 \To 1\) is not new, it follows e.g.~from the argument in the proof of \cite[Theorem 3.7]{enochs12}. For completeness, we still give a proof for \(2 \To 1\).

\begin{proof}[Proof of Theorem \ref{thm: duality of inj envelope and flat cover}]
\(1 \To 3\): Assume that \(i\) is an injective envelope. Let \(E\) be any injective module. Now \(\Hom_R(i,E)\) is an epimorphism, \(\Hom_R(I,E)\) is flat and cotorsion, and \(\ker \Hom_R(i,E) = \Hom_R(\coker i,E)\) is cotorsion. Thus, by Fact \ref{fact: flat cover iff top isom}, \(\Hom_R(i,E)\) is a flat cover, if and only if \(k(\mf p) \otimes_{R_\mf p} \Hom_R(R_\mf p,\Hom_R(i,E))\) is an isomorphism for all \(\mf p \in \spec R\). For any \(\mf p \in \spec R\) we have the commutative diagram
\[\begin{tikzcd}
k(\mf p) \otimes_{R_\mf p} \Hom_R(R_\mf p,\Hom_R(I,E)) \dar{k(\mf p) \otimes_{R_\mf p} \Hom_R(R_\mf p,\Hom_R(i,E))} \rar{\cong} & \Hom_R(\Hom_{R_\mf p}(k(\mf p),I_\mf p),E) \dar{\Hom_R(\Hom_{R_\mf p}(k(\mf p),i_\mf p),E)} \\
k(\mf p) \otimes_{R_\mf p} \Hom_R(R_\mf p,\Hom_R(M,E)) \rar{\cong} & \Hom_R(\Hom_{R_\mf p}(k(\mf p),M_\mf p),E)
\end{tikzcd}\]
where the horizontal morphisms are isomorphisms by Lemma \ref{lem: duality of soc and top}. Since the morphisms \(\Hom_{R_\mf p}(k(\mf p),i_\mf p)\) are isomorphisms, the morphisms \(\Hom_R(\Hom_{R_\mf p}(k(\mf p),i_\mf p),E)\) are isomorphisms as well. Hence the morphisms \(k(\mf p) \otimes_{R_\mf p} \Hom_R(R_\mf p,\Hom_R(i,E))\) are isomorphisms and thus \(\Hom_R(i,E)\) is a flat cover.

\(3 \To 2\): Trivial.

\(2 \To 1\): Assume that \(\Hom_R(i,E)\) is a flat cover for some injective cogenerator \(E\). By Fact \ref{fact: flat cover iff top isom}, the morphisms \(k(\mf p) \otimes_{R_\mf p} \Hom_R(R_\mf p,\Hom_R(i,E))\) are isomorphisms for all \(\mf p \in \spec R\). The above diagram shows that the morphisms \(\Hom_R(\Hom_{R_\mf p}(k(\mf p),i_\mf p),E)\) are all isomorphisms. Since \(E\) is an injective cogenerator, the morphisms \(\Hom_{R_\mf p}(k(\mf p),i_\mf p)\) are then all isomorphisms. Thus \(i\) is an injective envelope by Fact \ref{fact: injective envelope iff soc isom}.
\end{proof}

\begin{remark}
In \cite[Theorem 3.7]{enochs12} the authors give the following characterization of left Noetherian rings: a ring \(S\) is left Noetherian, if and only if a monomorphism \(i\) is an injective preenvelope, precisely when \(i^+\) is a flat precover. Here \((-)^+ := \Hom_\ZZ(-,\QQ/\ZZ) \cong \Hom_S(-,S^+)\); note that \(S^+\) is an injective cogenerator. If \(S\) is commutative, we can combine this characterization with Theorem \ref{thm: duality of inj envelope and flat cover} to see that \(S\) is Noetherian, if and only if a morphism \(i\) is an injective envelope, precisely when \(i^+\) is a flat cover.
\end{remark}

Theorem \ref{thm: duality of inj envelope and flat cover} of course immediately generalizes to minimal injective and flat resolutions.

\begin{corollary}\label{cor: duality of min inj reso and min flat reso}
Let \(0 \to M \to I^\bullet\) be a chain complex. The following conditions are equivalent:
\begin{enumerate}
\item \(0 \to M \to I^\bullet\) is a minimal injective resolution;
\item \(\Hom_R(I^\bullet,E) \to \Hom_R(M,E) \to 0\) is a minimal flat resolution for some injective cogenerator \(E\);
\item \(\Hom_R(I^\bullet,E) \to \Hom_R(M,E) \to 0\) is a minimal flat resolution for all injective modules \(E\).
\end{enumerate}
\end{corollary}

\begin{remark}
A referee pointed out that Corollary \ref{cor: duality of min inj reso and min flat reso} and thus Theorem \ref{thm: duality of inj envelope and flat cover} could alternatively be proved using the notion of minimality of complexes: a complex \(C\) is minimal, if and only if all homotopy equivalences \(C \to C\) are isomorphisms. We outline one alternative proof of the case \(1 \To 3\) using minimality.

First, note that if \(0 \to M \to I^\bullet\) is a minimal injective resolution, then the morphisms in the complex \(\Hom_{R_\mf p}(k(\mf p), I^\bullet_\mf p)\) are \(0\) for all \(\mf p\). Using Lemma \ref{lem: duality of soc and top} we get
\[k(\mf p) \otimes_{R_\mf p} \Hom_R(R_\mf p, \Hom_R(I^\bullet, E)) \cong \Hom_R(\Hom_{R_\mf p}(k(\mf p), I^\bullet_\mf p), E)\]
and see that the morphisms in the complex \(k(\mf p) \otimes_{R_\mf p} \Hom_R(R_\mf p, \Hom_R(I^\bullet, E))\) are thus \(0\) as well. Now, \cite[Theorem 2.3]{nakamura20} shows that \(\Hom_R(I^\bullet, E)\) is minimal. Finally, by \cite[Proposition 4.3(1)]{thompson19}, \(\Hom_R(I^\bullet,E) \to \Hom_R(M,E) \to 0\) is a minimal flat resolution.
\end{remark}

\begin{remark}
For an injective module \(E\) and a module \(M\), Corollary \ref{cor: duality of min inj reso and min flat reso} suggests a connection between the Bass numbers \(\mu_i(\mf p,M)\) and the dual Bass numbers \(\pi_i(\mf p,\Hom_R(M,E))\) defined by Xu in \cite[Definition 3.2]{xu95}. Indeed the following formula holds
\[\pi_i(\mf p,\Hom_R(M,E)) = \dim_{k(\mf p)} \Hom_R(k(\mf p),E)^{\mu_i(\mf p,M)},\]
where \(\Hom_R(k(\mf p),E)^{\mu_i(\mf p,M)}\) denotes the direct \emph{product} of copies of \(\Hom_R(k(\mf p),E)\) indexed by a set of cardinality \(\mu_i(\mf p,M)\). To see this, recall that the dual Bass numbers of the cotorsion module \(\Hom_R(M,E)\) satisfy
\[\pi_i(\mf p,\Hom_R(M,E)) = \dim_{k(\mf p)} \Tor_i^{R_\mf p}(k(\mf p),\Hom_R(R_\mf p,\Hom_R(M,E)))\]
by \cite[Theorem 2.2]{enochs97}. The formula follows then from the identity
\[\Tor_i^{R_\mf p}(k(\mf p),\Hom_R(R_\mf p,\Hom_R(M,E))) \cong \Hom_R(\Ext^i_{R_\mf p}(k(\mf p),M_\mf p),E).\]
This identity was already used in the proof of \cite[Proposition 2.9]{enochs97}, where the authors note that it implies \(\pi_i(\mf p,E(R/\mf p)) = \mu_i(\mf p,R)\).
\end{remark}

Next, we will consider the converse of Theorem \ref{thm: duality of inj envelope and flat cover}. It was proven in \cite[Corollary 3.3(2)]{mao15} that the converse holds with an added reflexivity requirement. For completeness, we give an easy proof of this using Theorem \ref{thm: duality of inj envelope and flat cover}.

\begin{proposition}[{\cite[Corollary 3.3(2)]{mao15}}]\label{prop: duality from flat to injective}
Let \(p \colon F \to M\) be a morphism where \(F\) is reflexive. Then \(p\) is a flat cover, if and only if \(p^\vee\) is an injective envelope.
\end{proposition}
\begin{proof}
Note that \(F\) is flat and \(p\) is an epimorphism, if and only if \(F^\vee\) is injective and \(p^\vee\) is a monomorphism. Thus we can assume that \(F\) is flat and \(p\) is an epimorphism.

Since \(F\) is reflexive and \(p\) is an epimorphism, \(M\) is reflexive as well. Hence, the commutative diagram
\[\begin{tikzcd}
F \rar{p} \dar{\cong} & M \dar{\cong} \\
(F^\vee)^\vee \rar{(p^\vee)^\vee} & (M^\vee)^\vee
\end{tikzcd}\]
shows that \(p\) is a flat cover, if and only if \((p^\vee)^\vee\) is a flat cover. The claim now follows from Theorem \ref{thm: duality of inj envelope and flat cover}.
\end{proof}

Both implications in the previous proposition fail if we drop the assumption that \(F\) is reflexive. We will show this in the following examples.

\begin{example}\label{ex: counterex 1}
Let \(f \colon F \to M\) be a flat precover. If \(\Hom_R(f,E)\) is an injective envelope for some injective cogenerator \(E\), then \(f\) is a flat cover. To see this, let \(h \colon F \to F\) be a morphism such that \(fh = f\). Now \(\Hom_R(h,E)\Hom_R(f,E) = \Hom_R(f,E)\) implies that \(\Hom_R(h,E)\) is an isomorphism. Thus \(h\) is an isomorphism.

This implication does not work if \(f\) is not a flat precover. For example, let \((R,\mf m)\) be a Noetherian local ring that is \emph{not} complete. We denote the residue field by \(k = R/\mf m\). Since \(\widehat{R} \to k\) is the flat cover of \(k\) and \(R \not\cong \widehat{R}\), \(f \colon R \to k\) is not a flat cover. However, \(f^\vee \colon k = k^\vee \to R^\vee = E(k)\) is an injective envelope.
\end{example}

\begin{example}\label{ex: counterex 2}
Let \((R,\mf m)\) be a Noetherian local ring that is not complete or \(\dim R \geq 2\). By \cite[Theorem 4.1]{belshoff92} there exists an injective envelope \(i \colon M \to I\) such that \(M\) is reflexive, but \(I\) is not. Now \(i^\vee\) is a flat cover by Theorem \ref{thm: duality of inj envelope and flat cover}, but \((i^\vee)^\vee\) is not an injective envelope since it factors through \(I \to (I^\vee)^\vee\) which is not an isomorphism.
\end{example}

\begin{example}\label{ex: counterex 3}
We can generalize the previous example further to non-local Noetherian rings. Let \(R\) be Noetherian with \(\mf m \in \maxspec R\) such that \(R_\mf m\) is not complete or \(\dim R_\mf m \geq 2\). Let \(f \colon F \to M\) be a flat cover in the category of \(R_\mf m\)-modules such that \(\Hom_{R_\mf m}(f,E(k(\mf m)))\) is not an injective envelope; this exists by the previous example. The morphism \(f\) is a flat cover in \(\RMod\) as well. To see this, note that any morphism \(F' \to M\) where \(F'\) is a flat \(R\)-module factors through \(F'_\mf m\). Thus there exists a morphism \(F'_\mf m \to F\) such that the diagram
\[\begin{tikzcd}
& F' \dar \\
& F'_\mf m \dar \ar{dl} \\
F \rar & M
\end{tikzcd}\]
commutes. Hence \(f\) is a flat precover in \(\RMod\). Clearly then \(f\) is a flat cover in \(\RMod\).

The morphism \(\Hom_{R_\mf m}(f,E(k(\mf m))) = \Hom_R(f,E(R/\mf m))\) is a direct summand of \(f^\vee\). Since \(\Hom_{R_\mf m}(f,E(k(\mf m)))\) is not an injective envelope, \(f^\vee\) is not an injective envelope either.

Even further, for \emph{any} injective cogenerator \(E\), the morphism \(\Hom_R(f,E)\) is not an injective envelope. This follows, since \(\bigoplus_{\mf n \in \maxspec R} E(R/\mf n)\) is a direct summand of \(E\) and so \(f^\vee\) a direct summand of \(\Hom_R(f,E)\).
\end{example}

We could of course directly extend Proposition \ref{prop: duality from flat to injective} to left resolutions that consist of reflexive modules. However, we can relax this assumption slightly with the following lemma.

\begin{lemma}
If \(M\) has a resolution consisting of reflexive flat modules, then the minimal flat resolution of \(M\) consists of reflexive modules as well.
\end{lemma}
\begin{proof}
Note that an epimorphism from a flat module with a cotorsion kernel is a flat precover. The kernels in a resolution consisting of reflexive flat modules are also reflexive and thus cotorsion. Hence, such a resolution is obtained by chaining flat precovers. Since the flat cover of any module \(N\) is a direct summand of any flat precover of \(N\) (see e.g.~\cite[Proposition 5.1.2]{enochs11}), the minimal flat resolution of \(M\) is a direct summand of any resolution of \(M\) obtained by chaining flat precovers. Therefore the minimal flat resolution also consists of reflexive modules.
\end{proof}

\begin{corollary}\label{cor: duality of reflexive min flat reso and min inj reso}
Assume that \(M\) has a resolution consisting of reflexive flat modules. Then, a chain complex \(F_\bullet \to M \to 0\) is a minimal flat resolution, if and only if \(0 \to M^\vee \to F_\bullet^\vee\) is a minimal injective resolution.
\end{corollary}
\begin{proof}
Assume first that \(F_\bullet \to M \to 0\) is the minimal flat resolution of \(M\). By the previous lemma, it consists of reflexive modules. Thus it is isomorphic to \((F_\bullet^\vee)^\vee \to (M^\vee)^\vee \to 0\). This being a minimal flat resolution implies that \(0 \to M^\vee \to F_\bullet^\vee\) is a minimal injective resolution by Corollary \ref{cor: duality of min inj reso and min flat reso}.

Assume then that \(0 \to M^\vee \to F_\bullet^\vee\) is the minimal injective resolution of \(M^\vee\). Again by Corollary \ref{cor: duality of min inj reso and min flat reso}, \((F_\bullet^\vee)^\vee \to M \to 0\) is a minimal flat resolution. The previous lemma then shows that \((F_\bullet^\vee)^\vee\) consists of reflexive modules. Since \(F_\bullet\) is a subcomplex of \((F_\bullet^\vee)^\vee\), \(F_\bullet\) also consists of reflexive modules. Thus \(F_\bullet = (F_\bullet^\vee)^\vee\) so \(F_\bullet \to M \to 0\) is a minimal flat resolution.
\end{proof}

\begin{remark}
The assumption of the previous corollary implies that \(M\) is reflexive. Note however that it is not enough to assume that \(M\) is reflexive. For example, extending the morphism \(i^\vee\) of Example \ref{ex: counterex 2} to a minimal flat resolution gives us a counterexample.
\end{remark}


\section*{Acknowledgements}

We thank the anonymous referees for their insightful comments. In particular for suggesting that the main result can also be proved using the notion of minimality of complexes. The author was supported by a grant from the Vilho, Yrjö and Kalle Väisälä Foundation.


\bibliographystyle{plain}
\bibliography{references}

\end{document}